\documentclass[10pt,a4paper]{article}
\newcommand{\MA}{\mathbb A}
\newcommand{\MD}{\mathbb D}
\newcommand{\Dom}{{\rm Dom}}
\newcommand{\zphinn}{\varphi_n}
\newcommand{\zphin}{\varphi_n(x)}

\newtheorem{Theorem}{Theorem}
\newtheorem{Corollary}[Theorem]{Corollary}
\newtheorem{Lemma}[Theorem]{Lemma}

\newtheorem{Remark}[Theorem]{Remark}
\newtheorem{Definition}[Theorem]{Definition}

\newcommand{\zaa}{\alpha}
 
\newcommand{\zg}{\gamma}
\newcommand{\ZDE}{\delta}
 \newcommand{\ZEP}{\epsilon}
\newcommand{\ZSI}{\sigma}
\newcommand{\zl}{\lambda} 
 
 \newcommand{\ZOMq}{\Omega} 
 \newcommand{\zthe}{\theta}
\newcommand{\zt}{\tau}

\def\zzn{\mathbb{N}}
\newcommand{\zzr}{\mathbb{R}}

\newcommand{\intt}{\int_0^t}

 \newcommand{\zdia}{~~\rule{1mm}{2mm}\par\medskip}
 
\newcommand{\ZIN}{\infty}
\newcommand{\zProof}{{\bf\underbar{Proof}.}\ }
 \newcommand{\ZD}{\;\mbox{\rm d}}
 \newcommand{\ZSUno}{\sum_{n=1}^{\ZIN}}

 \newcommand{\ZLA}{\label}

 \usepackage{epsfig}
\usepackage{graphics}
\usepackage[utf8]{inputenc}
\usepackage[english]{babel}
\usepackage{amsmath}
\usepackage{amsfonts}
\usepackage{amssymb}
\usepackage{makeidx}
\usepackage{graphicx}
\author{
L. Pandolfi\thanks{Retired from the Dipartimento di Scienze Matematiche ``Giuseppe Luigi Lagrange'', Politecnico di Torino, Corso Duca degli Abruzzi 24, 10129 Torino, Italy (luciano.pandolfi@formerfaculty.polito.it)}
}
\title{On the  fourth order Cattaneo equation of heat  conduction with memory}
\begin{document}
 
 \maketitle 
 
 \begin{abstract} 
The well known heat equation with finite speed of propagation proposed by Cattaneo is obtained by a more general fourth order PDE when a certain (small) parameters is put equal to zero. It seems that this fourth order equation has been essentially overlooked in the literature, in particular when it is subject to non homogeneous boundary conditions. In this paper we examine its well posedness   and the asymptotic behavior when certain coefficients tend to singular values.
\end{abstract} 
\section{Introduction}
An attempt to remediate the problem of infinite velocity of   diffusion of the standard heat equation led Cattaneo to propose its famous  equation  
\begin{equation}
\ZLA{eqCAtta}
\zt \zthe''=- a\zthe'+b\Delta\zthe\,,\qquad a>0\,,\quad b>0
\end{equation}
where $\zthe=\zthe(x,t)$ is the temperature at time $t$ and position $x\in\ZOMq \subseteq\zzr^d$   (in fact, Cattaneo considered the case $d=1$) and $\Delta$ denotes the laplacian, see~\cite{CattaneoMODENA1949}. This equation, which is a special instance of the telegraphists' equation, is now widely applied  in the description of diffusion processes in systems with complex molecular structure as seen from   a rapid search on the WEB (see for example~\cite{DEkEELIUHinestroza} for a description of the use of Cattaneo equation when modeling  diffusion of solutes in polymers)
and it is the prototype of the   heat equation with memory proposed by Gurtin and Pipkin in~\cite{GurtinPipkinARCHIVE1968} since it can be written as
\[
\zthe'(t)=b\int_{-\ZIN} ^t e^{-a(t-s) } \Delta \zthe(s)\ZD s\;.
\]

In the paper~\cite{CattaneoMODENA1949}, using statistical considerations, 
 Cattaneo   derived   a   general linear law  for the flux of heat:
\begin{equation}
\ZLA{eqFluxCatta}
q=-\chi\nabla  \zthe +\ZSI \nabla \zthe'  
\end{equation}
where $q=q(x,t)$ is  the flux, the prime denotes the time derivative and $\nabla$ is the gradient.
The parameter $\ZSI$ is small and if put equal to zero then we get the Fourier law which leads to the standard heat equation. Otherwise, as outlined in the next section, we get the fourth order equation
\begin{equation}
\ZLA{eqCAttaQUARTO}
 \zthe''=- a \zthe'+  b \Delta\zthe- c \Delta\zthe''  
 \qquad
 \left\{\begin{array}{l} 
  a= \chi/\ZSI >0\\
  b= \chi^2/(\ZSI\zg\rho)>0\\
  c= \ZSI/(\zg\rho )>0\,.
 \end{array}\right. 
\end{equation}

When the infinitesimal (for $\ZSI\to 0$) coefficient  $c$ is ignored then~(\ref{eqCAttaQUARTO}) reduces to the much studied equation~(\ref{eqCAtta}) while it seems that the complete equation~(\ref{eqCAttaQUARTO}) has not been considered, in particular in the case that the boundary conditions are nonhomogeneous.
Stimulated by the discussion in~\cite{SpiglerCATTANEOriv2020}, we
  present some observations on the \emph{fourth order} Cattaneo equation~(\ref{eqCAttaQUARTO}) when coupled with the following initial and boundary conditions:
\begin{equation}
\ZLA{eqCAttaQUARTOiniBOUND}
\zthe(0)=\zthe_0\,,\quad \zthe'(0)=\zthe_1\qquad \zg_0\zthe=f\ {\rm on}\ \Gamma=\partial\ZOMq\;.
\end{equation}
In fact, similar arguments can be adapted to study the case that the affine term $f$ acts in the Neumann boundary conditions.

We
 investigate whether the  problem~(\ref{eqCAttaQUARTO})-(\ref{eqCAttaQUARTOiniBOUND})   is
 well posed, i.e. whether it is  solvable, the solution is unique and depends continuously on the data in significant spaces. Roughly speaking, we prove:
\begin{itemize}
\item well posedness of Eq.~(\ref{eqCAttaQUARTO}) depends solely on the value of $c$. 
\item there exists an \emph{exceptional set} of values of the (positive) parameter c such that the problem \emph{is not} well posed (see Sect.~\ref{subS:NonWELl}). The exceptional set is the codomain of a sequence which tends to zero.
In terms of the parameter $\ZSI$, there exists a sequence $\{\ZSI_n\}$ which converges to zero and such that the system is not well posed if $\ZSI=\ZSI_n$.
\item otherwise, as proved in Sect.~\ref{sect:WELLposed}, the problem is well posed. In this case we find a representation for the solutions.
\item Similar to the standard heat equation, we prove that signals propagates with infinite speed (see Sect.~\ref{sect:finitePROPspeed}).

\item Finally, in Sect~\ref{secT:limits}, we investigate the behavior of the solutions when the parameters, either $c$ or $\ZSI$, tend to singular values. In particular, we study the case $\ZSI\to0^+$.
\end{itemize}
 
For completeness,  
in Sect.~\ref{Sect:DeriCatta} we present the arguments of Cattaneo without the deletion of any term, till to obtain Eq.~(\ref{eqCAttaQUARTO}).

\section{\ZLA{Sect:DeriCatta}The derivation of the  fourth order Cattaneo equation}
 
We outline the derivation of the Cattaneo equation as presented in~\cite{CattaneoMODENA1949}. Cattaneo first derives the  heat-flux law~(\ref{eqFluxCatta}) from statistical considerations and notes that the Fourier law is the special case obtained by putting $\ZSI=0$.   
  Putting $\ZSI=0$ removes ``memory'' from the heat processes and Cattaneo notes that this can be justified in a statical regime, but likely not before equilibrium is reached.
  To overcome this difficulty, instead of putting $\ZSI=0$ from the outset, Cattaneo couples Eq.~(\ref{eqFluxCatta}) with its time derivative, and obtains the system
\[
\left\{
\begin{array}{cl}
q&=-\chi\nabla  \zthe +\ZSI \nabla \zthe'\\
q'&=-\chi\nabla  \zthe'+\ZSI \nabla \zthe'' \;.
\end{array}
\right.
\] 
Elimination of $\nabla \zthe'$ gives
\begin{equation}
\ZLA{eqFluxCattaCOMPLE}
\ZSI q'=-\chi q-\chi^2\nabla\zthe+\ZSI^2\nabla\zthe''\;. 
\end{equation}
At this point Cattaneo notes that $\ZSI^2$ is infinitesimal of higher  order respect to $\ZSI $ and deletes the last term, proceeding with the computations below, which we perform without such semplification. Eq.~(\ref{eqFluxCattaCOMPLE})  is coupled with the conservation of energy:
\[
\left\{
\begin{array}{l}
\ZSI q' =-\chi q-\chi^2\nabla  \zthe+\ZSI^2 \nabla \zthe''  \\
\zg\rho \zthe'+\nabla\cdot q=0\;. 
\end{array}
\right.
\]
We equate the divergence of both the sides of the first equation and the time derivative of both the sides of the second equation. Using $\nabla \cdot q=-\zg\rho\zthe'$ (from the conservation of energy) we get the system
\[
\left\{
\begin{array}{l} 
\ZSI\nabla\cdot q' =-\chi\nabla\cdot q-\chi^2\Delta \zthe+\ZSI^2\Delta \zthe''=\chi\zg\rho\zthe'-\chi^2\Delta\zthe+\ZSI^2\Delta\zthe''\\
\zg\rho\zthe''+\nabla\cdot q'=0\;.
\end{array}
\right.
\]
Finally we replace $\nabla\cdot q'=-\zg\rho\zthe''$ from the second equation in the first one.  We get Eq.~(\ref{eqCAttaQUARTO}) which, we repeat, reduces to~(\ref{eqCAtta}) if, relaying on the ``smallness'' of $\ZSI$, the last addendum is removed.

\section{\ZLA{Sect:wellPOsedCattaFour}Well posedness and lack of well posedness of the fourth order Cattaneo equation}
We work in a bounded region $\ZOMq\subseteq \zzr^d$ with $C^2$ boundary $\Gamma$  and we introduce the \emph{Dirichlet laplacian}
 in $L^2(\ZOMq)$:
 \[
\Dom\,A=H^2(\ZOMq)\cap H^1_0(\ZOMq)\,,\qquad A\phi=\Delta\phi\quad \forall \phi\in\Dom\, A\;. 
 \]
It is known that $A$ 
  is selfadjoint negative defined with compact resolvent so that the space $L^2(\ZOMq)$ has an orthonormal basis $\{\zphinn\}$ of eigenvalues of $A$:
\[
A\zphinn=-\zl_n^2\zphinn\;.
\]  
Every eigenvalue $-\zl_n^2$ has finite multiplicity.

Let us introduce the \emph{exceptional set  for $c$}
\[
\mathcal{E}=\left \{\frac{1}{\zl_N^2}\,,\ -\zl_N^2\in\ZSI_p(A)\right \}.
\]
In this section we prove:
\begin{Theorem}\ZLA{teo:GeneBuonPosiz}
The following statements hold:
\begin{enumerate}
\item
\ZLA{I1teo:GeneBuonPosiz}
 if $c\in \mathcal{E}$ then problem~(\ref{eqCAttaQUARTO})-(\ref{eqCAttaQUARTOiniBOUND}) is not well posed in the sense that there exist initial conditions $\zthe_0$ and $\zthe_1$ such that no solution exist, not even with $f=0$.
\item
\ZLA{I2teo:GeneBuonPosiz} If $c\notin\mathcal{E}$ then problem~(\ref{eqCAttaQUARTO})-(\ref{eqCAttaQUARTOiniBOUND}) is well posed in the space $L^2(\ZOMq)\times L^2(\ZOMq)\times W^{2,p}(0,T;L^2(\Gamma))$ (any $p\geq 1$) in the sense that for every $\zthe_0$ and $\zthe_1$ in $L^2(\ZOMq)$ and  $f\in W^{2,p}(0,T;L^2(\Gamma))$, there exists a unique (mild) solution $\zthe\in C^1([0,T];L^2(\ZOMq))$ (any $T>0$) and the transformation $(\zthe_0,\zthe_1,f)\mapsto \zthe$  is continuous in the indicated spaces.
\end{enumerate}
 
\end{Theorem}

Note that the mild solutions  are defined in Sect.~\ref{sect:WELLposed}. 
   
   \subsection{\ZLA{subS:NonWELl}The proof of the statement~\ref{I1teo:GeneBuonPosiz} of Theorem~\ref{teo:GeneBuonPosiz}}

   We fix one of the eigenvalues $-\zl_N^2$ of $A$ and one of the corresponding eigenvectors, say $\varphi_N$:
   \[
A\varphi_N=-\zl_N^2\varphi_N\;.   
   \]
   
   Let 
   \[
c=\frac{1}{\zl_N^2}\;.   
   \]
   
  We impose the initial conditions
  \begin{equation}\ZLA{eq:Sec1:IniCondiPERnonWELL}
\zthe_0=\zaa \varphi_N\,,\qquad \zthe_1=\beta\varphi_N. 
  \end{equation} 
  The proof of of the statement~\ref{I1teo:GeneBuonPosiz} of Theorem~\ref{teo:GeneBuonPosiz} follows from the following observation:
  \begin{Lemma}
  Let $f=0$ and let $a>0$, $b>0$ be arbitrary while $c=1/\zl_N^2$.  The equation~(\ref{eqCAttaQUARTO}) with the initial conditions~(\ref{eq:Sec1:IniCondiPERnonWELL}) is solvable only if the coefficients $\zaa$ and $\beta$ satisfy the \emph{compatibility condition}\footnote{whose significance is clarified in Remark~\ref{rema:secTerzLIMIclaRIFcompat}.}
\begin{equation}
 \ZLA{eq:cOmpAtCond}
 \frac{\beta}{\zaa}=-\frac{b}{a}\zl_N^2\;.
 \end{equation}   
    \end{Lemma}
 \zProof  
  We proceed by separation of variables: a solution, if it exists, has the series expansion
  \[
\zthe(x,t)=\ZSUno \zphin\zthe_n(t)  
  \]
  and $\zthe_n(t)=0$ if $n\neq N$ while if $n=N$ then $\zthe_n(t)$ solves
  \[
  a\zthe_N'(t)=-b\zl_N^2\zthe_n(t )  \,.
  \]
  In fact, due to the equality $c=1/\zl_N^2$, the two terms containing $\zthe''$ cancel out.
  
 So that if the condition $\zthe_N(0)=\zaa$ has to be satisfied, it must be
  \[
  \zthe_N(t)=e^{-(b\zl_N^2/a) t}\zaa \;.
  \]
 The second condition $\zthe'(0)=\beta\varphi_N$ can be satisfied if and only if the    compatibility condition~(\ref{eq:cOmpAtCond})   holds.\zdia
  
   This ends the proof of the statement ~\ref{I1teo:GeneBuonPosiz} of Theorem~\ref{teo:GeneBuonPosiz}. This statement can be reformulated in terms of the physical parameters  of the system as follows. We introduce the \emph{exceptional set for $\ZSI$}
 \begin{equation}\ZLA{eq:defiEXCset}
 \mathcal{Z}=\left \{ \frac{\zg\rho}{\zl_n^2}\,,\qquad -\zl_n^2\in\ZSI_p(A)\right \}\;.
 \end{equation}
 Then we can state:
 \begin{Corollary}
 Let the coefficients $a$, $b$ and $c$ depend on $\ZSI$ as specified in~(\ref{eqCAttaQUARTO}).
 The problem~(\ref{eqCAttaQUARTO})-(\ref{eqCAttaQUARTOiniBOUND}) is not well posed if $\ZSI\in  \mathcal{Z}$.
 \end{Corollary}

 \subsection{\ZLA{sect:WELLposed}The proof of the statement~\ref{I2teo:GeneBuonPosiz} of Theorem~\ref{teo:GeneBuonPosiz}}
 In this section we prove that the system is well posed if $c\notin\mathcal{E}$.
As usual in PDE, the first step is to define the solutions in a suitable strong or mild sense. This can be easily achieved by noting that problem~(\ref{eqCAttaQUARTO})-(\ref{eqCAttaQUARTOiniBOUND}) can be rewritten as
\begin{equation}\ZLA{eq:CattaQuartPerBpoSIZ}
\left\{\begin{array}{l}
(1+c\Delta)\left (\zthe''(t)-(b/c)\zthe(t)\right )=h=-a\zthe'(t)-(b/c)\zthe(t)\,,\quad \mbox{in $\ZOMq$}\\
 \zg_0 \left (\zthe''(t)-(b/c)\zthe(t)\right )=f''-(b/c)f \quad \mbox{on  $\Gamma=\partial\ZOMq$}
 \end{array}\right.
\end{equation}
 and so, at least when $f\in W^{1,2}(0,T;L^2(\Gamma))$, for every fixed $t$ we have a Dirichlet problem which is well posed, since by assumption $-1/c$ is not an eigenvector of the Dirichlet laplacian.
 
 In order to proceed, we recall few known properties of this problem and we simplify the notations by introducing
 \[
H=L^2(\ZOMq)\,,\qquad U=L^2(\Gamma)   
 \]
and
\[
 \tilde A=I+cA\quad \mbox{on $H$ with $\Dom\,\tilde A=\Dom\,A =H^2(\ZOMq)\cap H^1_0(\ZOMq)$ (and $A=\Delta$)}\,.
\]
It is known:
\begin{equation}\ZLA{eq:soluLaplaEq}
(I+c\Delta)u=h\,,\qquad \zg_0 u=g\ \iff\ u=\tilde A^{-1}h+Dg
\end{equation}
where $D\in\mathcal{L}(U,H)$. 
Note that ${\rm im}\, D$ is not contained in $H^2(\ZOMq)$ and $u-Dg\notin\Dom\,\tilde A$ in general. 

By definition, $u=\tilde A^{-1}h+Dg$ is    a \emph{mild solution} of the Dirichlet problem~(\ref{eq:soluLaplaEq})   and the definition is justified by the fact that if $f$ is sufficiently smooth, and thanks to the assumption that $\partial\ZOMq$ is of class $C^2$, $u-Df\in \Dom\,\tilde A$ (see~\cite[Ch.IV par.~7]{mikailovLIBRO}).
We endow $\Dom\,\tilde A$ with the graph norm (which is equivalent to the $H^2(\ZOMq)$ norm
restricted to $\Dom\,\tilde A$).
It is known that   that $H $  has   continuous and dense injection in $(\Dom\,\tilde A)'$ and $\tilde A$ admits a continuous extension to $H$ (as a $(\Dom\,\tilde A)'$-valued function). Furthermore,  when $H=0$  and $f$ is smooth,
\[
(1+c\Delta)u=(1+c\Delta)(u-Df)=\tilde A u-\tilde ADf\qquad \mbox{in  $(\Dom\,\tilde A)'$} \,.
\]

Using these observations, when $f\in W^{2,2}(0,T;U)$ we can rewrite problem~(\ref{eq:CattaQuartPerBpoSIZ}) as
 \begin{equation}
\ZLA{eq:IIordOpeLimit}
\zthe''= (b/c) (I-\tilde A ^{-1})\zthe -a\tilde A^{-1}\zthe'- (b/c)  Df+Df''\;.
\end{equation}
This system can be put in the form of a standard semigroup system. We introduce the notations
\[
H=\zaa \tilde A^{-1}\,,\qquad K=\frac{b}{c}(I-\tilde A^{-1})\,,\qquad \beta=b/c
\]
and we have:
  \begin{equation}
  \ZLA{eq:PRIMOordOpeLimit}
   \frac{\ZD}{\ZD t}
 \underbrace{ \left(
  \begin{array}{c}
 \zthe \\ \zthe '
  \end{array}
  \right) }_{=W}
   =  \underbrace{ \left(
  \begin{array}{cc}
  0 &I\\
 K & -H
  \end{array}
  \right)}_{=\MA} W
 -\beta \underbrace{ \left(
  \begin{array}{c}
  0\\D
  \end{array}
  \right)}_{\mathbb{=D}} f
   +   \underbrace{
   \left(
   \begin{array}{c}
 0
\\D
   \end{array}
  \right)}_{=\MD} f''
=\MA W-\beta\MD f+\MD f''\,.
  \end{equation}

The form~(\ref{eq:CattaQuartPerBpoSIZ}) of the problem~(\ref{eqCAttaQUARTO})-(\ref{eqCAttaQUARTOiniBOUND})   and the previous observations 
suggest the following definition:
\begin{Definition}
\ZLA{DefiMILDsolu} 
A mild solution of the problem~(\ref{eqCAttaQUARTO})-(\ref{eqCAttaQUARTOiniBOUND}) is the first  component $\zthe$ of  a mild solution $W$ of the semigroup system~(\ref{eq:PRIMOordOpeLimit}). The solution is a classical solution when
the following properties hold :
\[
\begin{array}{l}
 \zthe(t)\in C^1([0,T];L^2(\ZOMq))\,,\quad   \zthe(0)=\zthe_0\,,\quad  \zthe'(0)=\zthe_1\,;\\
 y(t)=\zthe (t)-Df(t)\in C^2([0,T];L^2(\ZOMq))\,;\\
 y''(t) -(b/c) y(t)\in \Dom\,\tilde A\,;\\
 \tilde A(y''(t)-(b/c)y(t))= (b/c)\zthe(t)-a\zthe'(t)\;. 
\end{array}
\]
\end{Definition}
 
The study of the solution is now immediate, since
 $\MD
\in\mathcal{L}(U,H\times H)$ and
 $\MA\in\mathcal{L}(H\times H)$. So
 $\MA$   generates a \emph{uniformly continuous and holomorphic } group $e^{\MA t} $ of operators on $H\times H$ and we have
  \begin{equation}\ZLA{eqW(tPRIMdiIntegPARTI}
W(t)=e^{\MA t}W(0)  -\beta\intt e^{\MA(t-s)} \MD f(s)\ZD s+\intt e^{\MA(t-s)} \MD f''(s)\ZD s\;.
  \end{equation}
  The fact that $e^{\MA t}$ is a holomorphic   operator function shows that we can integrate by parts twice the second integral and we get:
  \begin{multline}\ZLA{eq:FormSEMIGRUPPI}
 W(t)= e^{\MA t}\left [
 W(0)-\MD f'(0)-\MA\MD f(0)
 \right ]\\
+\MD f'(t)+\MA\MD f(t)+\intt e^{\MA(t-s) } \left [-\beta I+\MA^2\right ]\MD f(s)\ZD s\;. 
  \end{multline}

We stress the following facts:
\begin{itemize}
\item the previous computations used explicitly $c\notin\mathcal{E}$;
\item the computations used $f\in W^{2,p}(0,T;L^2(\Gamma))$ (any $p\geq 1$ and any $T>0$)
but the right hand side makes sense under the weaker assumption that $f'$ is continuous. So, the mild solutions are defined under this weaker condition.
\end{itemize}
So:
 \begin{Theorem}\ZLA{teo:stateSpaziWellPos}
 Let $c\notin\mathcal{E}$. Then the equation~(\ref{eqCAttaQUARTO}) is well posed in the space $H\times H\times \in W^{2,p}(0,T;L^2(\Gamma))$ in the sense 
  that for every $\zthe_0$ and $\zthe_1$ in $H=L^2(\ZOMq)$ and every $f \in W^{2,p}(0,T;L^2(\Gamma))$
 there exists a unique mild solution $\zthe\in C^1([0,T];H)$ and the transformation
 \[
(\zthe_0,\zthe_1,f)\mapsto \zthe 
 \] 
 is linear and continuous from $H\times H\times W^{2,p}(0,T;L^2(\Gamma))$ to $C^1([0,T];H)$ for every $T>0$ and also from $H\times H\times C^1([0,T];L^2(\Gamma))$ to $C^1([0,T];H)$
 (where $H=L^2(\ZOMq)$).
   \end{Theorem}

 In terms of the physical parameters, we have:
 
 \begin{Corollary}
 Eq.~(\ref{eqCAttaQUARTO}) is well posed in the spaces specified in Theorem~\ref{teo:stateSpaziWellPos} if and only if $\ZSI\notin\mathcal{Z}$ (the   set $\mathcal Z$ is defined in~(\ref{eq:defiEXCset})). 
 \end{Corollary}

Finally, we justify in the usual manner the definition of the mild solutions:
\begin{Theorem}
Any mild solution  $\zthe(t)$   is the $C^1([0,T];H)$-limit of a sequence of classical solutions.
\end{Theorem} 
\zProof
 It is sufficient to prove that the solution is classical when $(\zthe_0,\zthe_1,f)$ belong to a dense subset of $L^2(\ZOMq)\times L^2(\ZOMq)\times  W^{2,p}(0,T;L^2(\Gamma))$. Using the linearity of the problem, we can study separately the initial condition and the boundary function $f$. We confine ourselves to examine the contribution of $f$ leaving the contribution of the initial conditions to the reader.
 
 We assume $f\in C^\ZIN([0,T]\times\Gamma)$ and we consider equality~(\ref{eqW(tPRIMdiIntegPARTI}) (with $W(0)=0$). The regularity of $f(t)$ shows that
 $W(t)\in C^1([0,T],L^2(\ZOMq)\times L^2(\ZOMq))$ (in fact also in $C^k$ for every $k$)
 and the following equality holds for every $t$:
 \[
W'(t)=\MA W(t)-\beta \MD f(t)+\MD f''(t) 
 \]
so that
\[
W'(t)+\beta \MD f(t)-\MD f''(t)\in {\rm im}\, \MA 
\]
and the equality~(\ref{eq:PRIMOordOpeLimit}) holds for $t\in[0,T]$.
The second component of the equality~(\ref{eq:PRIMOordOpeLimit}) gives
\[
\zthe''(t)-Df''-\beta(\zthe-Df)=\tilde A^{-1}(-\zaa\zthe+\beta\zthe)\in\Dom\tilde A \;.
\]
From this equality, the properties required    in Definition~\ref{DefiMILDsolu} are easily seen.\zdia
 
\begin{Remark}{\rm 
We note the similarity of Eq.~(\ref{eqCAttaQUARTO})  (when it is well posed) with the standard heat equation: when $f=0$ the solutions are analytic functions of $t$ while the similarity with the wave equation is that the equation is reversible and the smoothness of the initial condition $W(0)$ is preserved.  }
\end{Remark} 
 
The following observation
is similar to that used in~\cite{bucci-las-optimiz_2018,LasiePandLukes2,LasiePandTrigg}. 
It might have an interest, although its physical significance is not  clear in the contest of Eq.~(\ref{eqCAttaQUARTO}).  In the cited papers, this observation was used 
to solve a quadratic optimization problem.

When $f\in W^{2,p}(0,T;L^2(\Gamma)) $, the formula~(\ref{eq:FormSEMIGRUPPI}) decouples $f(0)$ and $f'(0)$ from $f(t)$, $t>0$ and introduces the parameters $f_0$ (i.e. $f(0)$), $f_1$ (i.e. $f'(0)$).  After this decoupling, formula~(\ref{eq:FormSEMIGRUPPI}) makes sense if $f'$  is solely (square) integrable.   Once this decoupling has been done, $W(t)$ is not continuous. If $f\in W^1(0,T;L^2(\Gamma))$ then the formula may be used to define a solution $W(t)$ in an even weaker sense: it not continuous: as an element of $L^2(0,T,H\times H)$ depends continuously
on f $f\in W^1(0,T;L^2(\Gamma))$ and on
the parameter  $W_1=W(0)-\MD f_1-\MA\MD f_0$. Its first component $\zthe(t)$ is given by the equality
\begin{multline*}
\zthe(t)=\left [\begin{array}{ll} I &0\end{array}\right] \MA^{-1}\left (e^{\MA t}-I\right )W_1+Df(t)\\
+\intt
\left [\begin{array}{ll} I &0 \end{array}\right] e^{\MA(t-s)} [-\beta I+\MA^{2}]\MD f(s)\ZD s
\end{multline*}
where $0$ and $I$ denote the zero and identity operators in $H$. This solution $\zthe\in L^2([0,T];H)$ depends continuously on $W_1$ and on $f\in L^2(0,T;U)$ and furthermore $\zthe(t)-Df(t)$ is continuous.

\section{\ZLA{sect:finitePROPspeed}Propagation speed}
 One of the reason for the introduction of Eq.~(\ref{eqCAtta}) was that waves propagates with finite speed in a medium described by Eq.~(\ref{eqCAtta}). Instead, we are going to show the 
 existence of  solutions of Eq.~(\ref{eqCAttaQUARTO})   whose propagation speed is as great as we wish \emph{when the coefficient $c$ is strictly positive and $c\notin\mathcal{E}$}.
Note that this fact is easily guessed in the case that $\ZOMq=\zzr$. In this case  for every fixed $t>0$ we can formally compute the Fourier transform $\hat \zthe(\zl,t)$ of the solution. We have
\begin{equation}\ZLA{ODEdopoTransf}
(1-c\zl^2)\hat \zthe''(\zl,t) +a\hat \zthe'(\zl,t) +\zl^2b\hat \zthe(\zl,t)=0\qquad \left\{\begin{array}{l}\hat\zthe(\zl,0)=\hat w_0(\zl) \\ \hat\zthe'(\zl,0)=\hat w_1(\zl) \end{array}\right.
\end{equation} 
(the apex denotes the time derivative).

Let $w_0=0$ and let $w_1$ have compact support. If signals propagate with finite speed then 
  $x\mapsto w(x,t)$ has compact support for every $t>0$. Then the function $\zl\mapsto \hat w(\zl,t)$ has to be entire (and in fact it must have stronger properties, see Paley-Wiener Theorem in~\cite[p.~132]{KoosisLIBROhP}).

Eq.~(\ref{ODEdopoTransf}) is an ordinary differential equation for every fixed value of $\zl$, whose solution is
 \begin{multline}\ZLA{eq:soluConTrasFOURIE}
\hat\zthe(\zl,t)=A(\zl)\exp\left \{ \frac{-a+\sqrt{a^2-4b\zl^2(1-c\zl^2)}}{2(1-c\zl^2)}t \right  \}\\-
A(\zl)\exp\left \{ \frac{-a-\sqrt{a^2-4b\zl^2(1-c\zl^2)}}{2(1-c\zl^2)}t \right \} 
 \end{multline}
 and
 \[
A(\zl)=\frac{1-c\zl^2}{\sqrt{a^2-4b\zl^2(1-c\zl^2)}} \hat w_1(\zl)\;.
 \]
The function $\hat w_1(\zl)$ is entire and so the function $A(\zl)$ has a zero of finite order  for $\zl=\pm \sqrt{c}$. It cannot compensate the essential singularity for $\zl=\pm\sqrt c$ of the second exponential in~(\ref{eq:soluConTrasFOURIE}) (the first exponential remains bounded  for $\zl\to \pm\sqrt c$ as it is easily seen by expanding the square root with   the binomial formula).

 An argument which applies also to bounded regions 
   is easily seen from formula~(\ref{eq:FormSEMIGRUPPI}). Let $W_0=0$. We fix any $T>0$ and we  consider the sequence of $C^{\ZIN}$ input functions
\[
f_n(t)=\frac{1}{n}e^{-n (T-t)}  f_0\,,\qquad \mbox{any $f_0\in U=L^2(\Gamma) $ such that $\MD f_0\neq 0$}\,.
\]
Note that
\[
\begin{array}
{l}
\displaystyle    \lim _{n\to+\ZIN} f_n(t)=0\quad \mbox{uniformly for $t\in[0,T]$}      \\
\displaystyle\lim _{n\to+\ZIN} f_n'(t)=0 \quad \mbox{in  $L^2(0,T;U)$ (uniformly for $t\in[0,T-\ZEP]$, any $\ZEP>0$)}\\
\displaystyle f_n'(T)= f_0\quad \mbox{for every $n$}\,.
\end{array}
\]
Let $W_n(t)$ be the solution. Formula~(\ref{eq:FormSEMIGRUPPI}) gives
\begin{equation}\ZLA{LimiPerVelOInfi}
\lim _{n\to+\ZIN} W_n(T)= \MD f_0\quad \mbox{in $H=L^2(\ZOMq)$}\;.
\end{equation}
The element $\MD f_0\in H$  is non zero. Hence there exists a (nonempty) subregion $\ZOMq_1\subseteq \ZOMq$ (so $\ZOMq_1$ of positive measure) \emph{at positive distance $d$ from $\partial\ZOMq$} such that
\[
\int _{\ZOMq_1} \|\MD f\|^2\ZD x=\zaa>0\;.
\]
  The property~(\ref{LimiPerVelOInfi}) shows the existence of  $ N(T)$ such that if $n>N(T)$ then we have  
\[
\int _{\ZOMq_1} \|w_n(T)\|^2\ZD x=\zaa/2>0  \;.
\]
We recall that we can fix  $T>0$, as small as we wish.
So, for every (small) $T>0$ there exist waves excited on the boundary   which cover the distance $d$ in time less then $T$.

\section{\ZLA{secT:limits}The solutions when $c$ and $\ZSI$ tends to exceptional values}

Eq.~(\ref{eqCAtta}) is the special case of Eq.~(\ref{eqCAttaQUARTO}) when the coefficient $c$ is put equal zero thanks to the fact that it is small when compared with the other coefficients. It is not 
stated in~\cite{CattaneoMODENA1949} that the case $c=0$ is obtained as the limit of~(\ref{eqCAttaQUARTO}) for $c\to 0$. And in fact, we are going to see that in general the  solutions, as a function of $c$, do not have limits. Even more so because the parameter which should tend to zero is the physical parameter $\ZSI$. In terms of $\ZSI$, Eq.~(\ref{eqCAttaQUARTO}) is the equation
\begin{equation}
\ZLA{eqCAttaQUARTOZSI} 
\ZSI\zthe''+a\zthe'=b\Delta\zthe-\ZSI c\Delta \zthe'' \quad \mbox{where we redefined} \quad \left\{
\begin{array}{l}
a=\chi\\
b=\chi^2/(\zg\rho)\\
c=1/(\zg\rho)\;.
\end{array}
\right.
\end{equation}
From this point of view, we might expect that the solutions of Eq.~(\ref{eqCAttaQUARTO}) when $\ZSI\to 0^+$ converge to those of the heat equation
\[
a\zthe'=b\Delta\zthe                                                       
\]
at least when the initial conditions satisfy the compatibility  condition~(\ref{eq:cOmpAtCond}). We disprove also this conjecture.

We sum up: in this section we consider three significant limit processes:
\begin{enumerate}
\item when $a$, $b$ are kept fixed and $c\to 1/\zl_n^2$;
\item when $a$ and $b$ are kept fixed and  $c\to 0$ while being $c\notin\mathcal{E}$;
\item when the complete dependence on $\ZSI$ is considered and so the equation is~(\ref{eqCAttaQUARTOZSI}). We prove the existence of solutions  which do not converge when $\ZSI\to 0^+$ with the 
condition $\ZSI\notin\mathcal{Z}$ (the set $\mathcal{Z}$ is defined in~(\ref{eq:defiEXCset})) not even if the initial conditions $w_0$ and $w_1$ satisfy the compatibility condition
$aw_1=b\Delta w_0$
 imposed by the heat equation.
\end{enumerate}

The interesting case are the second and third ones but the computation in the first case provides preliminary formulas.

\subsection{The first limiting process: when $c\to 1/\zl_n^2$ with $c\notin\mathcal{E}$}
We recall the definition of the \emph{exceptional set}
\[
\mathcal{E}=\left \{\frac{1}{\zl_n^2}\,,\ -\zl_n^2\in\ZSI_p(A)\right \}
\]
and we recall that the problem is well posed if and only if $c\notin\mathcal{E}$. We show the existence of initial conditions such that the corresponding solution does not have limit when $c\to 1/\zl_n^2$ for a fixed value of $n$. An obvious guess is that the limit should not exist if the initial condition does not satisfy the compatibility condition~(\ref{eq:cOmpAtCond}), for example if
 
 \[
f=0\,,\qquad \zthe_0=0\,,\quad \zthe_1=\zphinn 
 \]
 where $\zphinn$ is an eigenvector of the eigenvalue $\zl_n^2$. This is true and we leave the verification to the reader since we can prove a stronger and possibly unexpected property: the limit does not exist even if the compatibility condition is satisfied. So, we consider the solution with $f=0$ and
 \[
\zthe(0)=-\frac{a}{b\zl_n^2}\zphinn\,,\qquad \zthe'(0)=\zphinn\;. 
 \]
 
 Proceeding by separation of variables as in Sect~\ref{subS:NonWELl} but with $c\neq 1/\zl_n^2$  we see that the solution is
 \[
\zthe(t)=\zthe_n(t)\zphinn =\underbrace{A\exp\left \{\frac{-a+\ZDE }{2(1-c\zl_n^2)}t\right \}}_{\zthe_{1,c}(t)}
 +\underbrace{B\exp\left \{-\frac{a+\ZDE }{2(1-c\zl_n^2)}t\right \}}_{\zthe_{2,c}(t)}
 \]
 where
 \begin{multline}\ZLA{defiDELradicandoZDE}
 \ZDE =\sqrt{a^2-4b\zl_n^2(1-c\zl_n^2)}\\
 =a\left [
 1-\frac{2b\zl_n^2}{a^2}(1-c\zl_n^2)-\frac{2b^2\zl_n^4}{a^4}(1-c\zl_n^2)^2+{\rm o}(1-c\zl_n^2)^2
 \right ]\;.
 \end{multline}
 The coefficients $A$ and $B$ are found by imposing the initial condition to $\zthe_n(t)$:
\[ 
\zthe_n(0)=-a/b\zl_n^2\,,\qquad \zthe_n'(0) =1 
\]
 and it turns out that
 \[
A=\frac{2b\zl_n^2(1-c\zl_n)^2-a(a+\ZDE)}{2b\ZDE\zl_n^2}\,,\qquad B=-\frac{a(\ZDE-a)+2b\zl_n^2(1-c\zl_n^2)}{2b\ZDE\zl_n^2}\;.
 \]
In these expressions, $n$ and $\zl_n^2$ are fixed while we consider the limit for $c\to 1/\zl_n^2$.
It is easily seen that
\[
\left\{\begin{array}{l}
\lim _{c\to 1/\zl_n^2} A=-a/b\zl_n^2 \\
\lim _{c\to 1/\zl_n^2} B=0\\
{\rm and}\ B=\frac{b\zl_n^2}{a^3}(1-c\zl_n^2)^2+{\rm o}(1-c\zl_n^2) 
\end{array}\right.
\quad 
 \left\{\begin{array}{l}
 \lim  _{c\to (1/\zl_n^2)} \frac{-a+\ZDE}{2(1-c\zl_n^2)}=-2b\zl_n^2/a \\[2mm]
\lim _{c\to (1/\zl_n^2)^+}   [-\frac{a+\ZDE }{2(1-c\zl_n^2)}  ]=+\ZIN 
\\[2mm]
\lim _{c\to (1/\zl_n^2)^-}  [-\frac{a+\ZDE }{2(1-c\zl_n^2)}  ]=-\ZIN\,.
\end{array}\right.
 \]
So we have for $t>0$:
 \[
\left\{\begin{array} 
{l}
\lim _{c\to 1/\zl_n^2} \zthe_{1,c}(t)=-\frac{a}{b\zl_n^2}e^{-(2b\zl_n^2/a)t}\\
\lim _{c\to (1/\zl_n^2)^+} \zthe_{2,c}(t)=\ZIN \\
\lim _{c\to (1/\zl_n^2)^-} \zthe_{2,c}(t)=0\,.
 \end{array}\right.
 \]
  
\subsection{The second limiting process: when $c\to 0$  with the condition  $c\notin\mathcal{E}$}

Here we show the existence of solutions of Eq.~(\ref{eqCAttaQUARTO}) which do not converge when $c\to 0^+$. An example is as follows:

Let the initial condition be 
\[
\zthe_0=0\,,\qquad \zthe_1=\ZSUno \frac{1}{n}\zphinn\;.
\]
Then, with $\zthe_{n,c}$ the solution in~(\ref{Eq:SuluZThecN}),  the solution $\zthe=\zthe_c(t)$ is
\[
\zthe_c(t)=\ZSUno\frac{1}{n} \zphinn\zthe_{n,c }(t)  
\]
where $\zthe_n=\zthe_{n,c}(t)$  solves
\begin{equation}
\ZLA{Eq:SuluZThecN}
(1-c\zl_n^2)\zthe_n''+a\zthe_n'+b\zl_n^2\zthe_n=0\qquad \zthe_n(0)=0\,,\quad \zthe_n'(0)=1\;.
\end{equation}
The solution is
\begin{equation}
\ZLA{Eq:SuluZThecNsuaSOLU}
\zthe_{n,c}(t)=\frac{ 1-c\zl_n^2}{ \ZDE}\left [\exp\left \{\frac{-a+\ZDE }{2(1-c\zl_n^2)}t\right \}-
\exp\left \{-\frac{a+\ZDE }{2(1-c\zl_n^2)}t\right \}\right ]
\end{equation}
and $\ZDE$ is in~(\ref{defiDELradicandoZDE}).

We want to study the limit of $\zthe_c(t)$ when $c\to 0$ while staying outside of the exceptional set $\mathcal{E}$.

For every $t$ and every fixed $n$ we have
\begin{equation}\ZLA{eq:sectLimIProceDisegSUP}
 \|\zthe_{ c}(t)\|_{H}\geq \frac{1}{n}|\zthe_{n,c}(t)|\quad \mbox{so that}\quad  \|\zthe_{c}(t)\|_{H}\geq\sup\left \{ \frac{1}{n}|\zthe_{n,c}(t)|\right \}\;.
\end{equation}
We prove
\begin{equation}\ZLA{eq:supDAprovSecoLim}
\lim _{c}\sup\left \{ \frac{1}{n}|\zthe_{n,c}(t)|\right \} =+\ZIN
\end{equation}
the limit being computed while $c\to 0$ along the sequence
\[
c_k=\frac{1}{	\zl_k^2}+\zg\frac{1}{\zl_k^3}\;.
\]
The number $\zg>0$ is fixed and chosen in such a way that $c_k\notin\mathcal{E}$. The number  $\zg$ exists since $\mathcal{E}$ is denumerable.

In order to verify~(\ref{eq:supDAprovSecoLim}) it is sufficient that we prove
\[
\lim _{c=c_k\to 0^+}\frac{1}{k}|\zthe_{k,c_k}(t)|=+\ZIN\,.
\]

We note that
\[
c_k\zl_k^2=1+\frac{\zg}{\zl_k},,\quad 1-c_k\zl_k^2=-\frac{\zg}{\zl_k}\,,\qquad \ZDE=\sqrt{a^2+4b
\zg\zl_k}\;.
\]

The common coefficient of $\zthe_{k,c_k}(t)$ is
\[
\frac{1}{k}\frac{1-c\zl_k^2}{\ZDE}=-\frac{\zg}{k\zl_k}\frac{1}{\sqrt{a^2+4b\zg\zl_k}}\to 0\;.
\] 
Weyl estimate of the eigenvalues of the Dirichlet laplacian in a bounded region of $\zzr^d$ (see~\cite{62CourantHilbertBOOK2,mikailovLIBRO}) gives
\[
\zl_k\asymp k^{1/d}\qquad  d={\rm dim}\,\ZOMq
\]
so that
\begin{equation}\ZLA{eq:secoLIMasimpCOeff}
\frac{1}{k}\frac{1-c\zl_k^2}{\ZDE}\asymp\frac{1}{k^{1+( 3 /2d)}}\;.
\end{equation}

The first exponential remains bounded since the exponent is negative:
\[
\frac{-a+\ZDE}{2(1-c\zl_n^2)}t=\frac{-a+\sqrt{a^2+4b\zg\zl_k}}{-2\zg}\zl_kt\leq 0\;.
\]
So, the contribution of the first addendum tends to zero. Instead, we prove that the second addendum diverges for every fixed $t>0$. Thanks to the estimate~(\ref{eq:secoLIMasimpCOeff}), it is sufficient that we examine the exponent, which is
\[
-\frac{a+\ZDE}{2(1-c\zl_n^2)}t=-a\frac{1+\sqrt{1+4b\zg\zl_k/a^2}}{-2\zg }\zl_kt\asymp k^\frac{3}{2d}\;.
\]
The required property~(\ref{eq:supDAprovSecoLim}) follows.

This observation confirms the fact that Eq.~(\ref{eqCAtta}) is not the limit of~(\ref{eqCAttaQUARTO}).

\subsection{\ZLA{sec:TerzoLimITE}The third limiting process: when $\ZSI\to 0$}
 
 In spite of the fact that this may look as the most natural limit process to be considered,  also in this case we prove the existence of solutions which do not converge.
 
 We consider the following example:
 \[
\ZOMq=(0,\pi)  \quad \mbox{so that $\zl_n^2=-n^2$ (simple eigenvalue)}
 \]  
 and the problem
\begin{equation}
 \ZLA{eq:step3delLIMIT}
\ZSI \zthe''=-2\zthe'+\Delta\zthe-\frac{\ZSI^2}{4}\Delta\zthe''\,,\qquad \left\{\begin{array}
{l}
\displaystyle 
\zthe(0)=\ZSUno \frac{1}{n^4}\zphinn\,,\\
\displaystyle 
\zthe'(0)=-\ZSUno \frac{1}{2n^2} \zphinn \,,\\
\displaystyle 
\zthe(0,t)=\zthe(\pi,t)=0 
\end{array}\right.
 \end{equation}
 
\begin{Remark}\ZLA{rema:secTerzLIMIclaRIFcompat}
{\rm We note:
\begin{itemize}
\item
  the initial conditions $\zthe_0$ and $\zthe_1$ satisfy the compatibility condition imposed by the heat equation $2\zthe'=\Delta\zthe$:
 \[
2\zthe_1= 2\zthe'(0)=\Delta\zthe(0)=\Delta\zthe_0
 \]
 and this compatibility condition imposed by the heat equation is precisely the condition that the components $\zthe_{0,n} $ and $\zthe_{1,n}$ of the initial data satisfy the compatibility conditions~(\ref{eq:cOmpAtCond}).
\item   the coefficients of the system are obtained when the values of the physical parameters are $\chi=2$ and $\zg\rho=4$,
compare~(\ref{eqCAttaQUARTO}) and~(\ref{eqCAttaQUARTOZSI}). So, this example is not artificial.
\end{itemize}
}
\end{Remark}

The condition $\ZSI\to 0$  and $\ZSI\notin\mathcal{Z}$ is the condition
\begin{equation}
\ZLA{eq:UltimoLimCINDIZSI}
\ZSI\to 0\,,\qquad \ZSI\neq\frac{4}{k^2}\quad\forall k\;.
\end{equation}
We shall choose the sequence $\{\ZSI_k\}$
\[
\ZSI_k=\frac{5}{k^2}   
\]
 so that $ \ZSI\notin\mathcal{Z}$ since $\sqrt 5$ is irrational.

We solve~(\ref{eq:step3delLIMIT}) by separation of variables:
\[
\zthe(x,t)=\ZSUno\zphin \zthe_n(t)
\] 
where $\zthe_n(t)$ solves
\[
\left (\ZSI-\frac{1}{4}\ZSI^2n^2\right )\zthe''=-2\zthe_n'-n^2\zthe_n\qquad \zthe_n(0)=\frac{1}{n^4}\,,\quad \zthe_n'(0)=-\frac{1}{2n^2}\;.
\]
 
It is easily computed that
\[
\zthe_n(t)=-\frac{1}{8n^4}\frac{(4-\ZSI n^2)^2}{n^2\ZSI-2}\exp\left\{ -\frac{2n^2}{4-\ZSI n^2}t\right\}+\frac{ \ZSI^2}{8( n^2\ZSI-2)} \exp\left \{-\frac{2}{\ZSI}t \right \}\;.
\]
 
Now we use a similar idea to that used in the study of the second limit. We fix any $t>0$ and we denote explicitly the dependense of $\ZSI$:
\[
\zthe(x,t)=\zthe_\ZSI(x,t)\,,\quad \zthe_n(t)=\zthe_{n,\ZSI}(t)\,.
\]
Then, for every $\ZSI>0$, \  $\ZSI\notin\mathcal{Z}$, we have
\begin{equation}\ZLA{secIIIlimiteEq}
 \|\zthe_\ZSI(t,\cdot)\|_{L^2(0,\pi)}> |\zthe_{n,\ZSI}(t)|\quad \mbox{and so}\quad 
  \|\zthe _\ZSI(t,\cdot)\|_{L^2(0,\pi)}>\sup \{  |\zthe_{n,\ZSI}(t) |\,,\quad n\in\zzn\}\;.
\end{equation}
Now we evaluate both the addenda of $  \zthe_{k,\ZSI}(t)  $ when $\ZSI=\ZSI_k=5/k^2$.
\begin{itemize}
\item the first addendum is $(1/24k^2)e^{2k^2 t}\to+\ZIN$. We have  $(1/24k^2)e^{2k^2 t}> 2k$ if $k$ is large.
\item the second addendum is $(25/24 k^4)e^{-2k^2t/5}<2/k^4<2$.
\end{itemize}
And then, there exists $k$ such that  
\[
|\zthe_{k,\ZSI_k}(t)|_{L^2(0,\pi)}>k\qquad \mbox{so that}\quad  \|\zthe_{5/k^2}(t,\cdot)\|_{L^2(0,\pi)}>\to+\ZIN\;.   
\]
It follows that  $\zthe_{5/k^2}(t)$ do not converge to a solution of the standard heat equation.

    \bibliography{bibliomemoria}{ }
  \bibliographystyle{plain}
 \enddocument